\documentclass[11pt]{amsart}
\usepackage{amsmath}
\usepackage[T1]{fontenc}
\usepackage[utf8]{inputenc}
\usepackage{graphicx} 
\usepackage{subfig}   
\usepackage{geometry}                
\usepackage{amssymb}
\usepackage{amsthm}
\usepackage{enumerate}
\usepackage{hyperref}
\usepackage{pinlabel}
\hypersetup{
urlcolor= blue, 
linkcolor= blue, 
bookmarksopen=true, 
pdftitle={Transitive parallelism of residues in buildings}, 
pdfauthor={Antoine Clais}, 
pdfsubject={Mac OS X} 
}

\newtheorem{theo}{Theorem}
\numberwithin{theo}{section}
\newtheorem{prop}[theo]{Proposition} 
\newtheorem{coro}[theo]{Corollary}
\newtheorem{lem}[theo]{Lemma} 
\newtheorem{déf}[theo]{Definition}

\newtheorem{exmp}[theo]{Example} 

\newcommand{\Z} {\ensuremath{\mathbb{Z}}}
\newcommand{\h} {\ensuremath{\mathbb{H}}}
\newcommand{\N} {\ensuremath{\mathbb{N}}}

\newcommand{\s} {\ensuremath{\mathbb{S}}}

\newcommand{\proj}[2] {\ensuremath{\mathrm{proj}_{#1} (#2)}}
\newcommand{\apD} {\ensuremath{\mathcal{A}p(\Delta)}}

\newcommand{\cox} {\ensuremath{(W,S)}}
\newcommand{\coxsyst} {\ensuremath{\Sigma(W,S)}}
  
 \newcommand{\Autt}[1] {\mathrm{Aut}_{\mathrm{T}}(\ensuremath{#1})}   
\newcommand{\di}[1] {\mathrm{dist}(\ensuremath{#1})}

\newcommand{\dc}[1] {\ensuremath{d_c}(\ensuremath{#1})}
\newcommand{\dia}[1] {\mathrm{diam} \ensuremath{\:#1}}

   \newcommand{\numeroti} [1] {
  \begin {enumerate} [i)]
  #1
  \end{enumerate}}
  
     \newcommand{\liste} [1] {
  \begin {itemize}#1
  \end{itemize}}

\title{Transitive parallelism of residues in buildings}
\author{Antoine Clais}
\address{Laboratoire Paul Painlevé \\ Université Lille 1 \\ 59655 Villeneuve d'Ascq, France} 
\email{antoine.clais@math.univ-lille1.fr}
\date{\today} 

\begin{document}

\begin{abstract}
We study the buildings in which parallelism of residues is an equivalence relation. If the building admits a group action, we  describe how parallel residues are related to residues with equal stabilizers.  This permits to retrieve the fact that in a Coxeter group or in a graph product, intersections of parabolic subgroups are parabolic.
\end{abstract}

  \maketitle
 \begin{flushleft}{\bf Keywords:} Buildings, Coxeter groups, parallel residues, parabolic subgroups.
 
\end{flushleft}
 \section{Introduction}

 In a building $\Delta$, residues are convex subsets equipped with a natural building structure directly inherited from $\Delta$. In \cite{TitsBuildingsLectureNotes}, J. Tits has introduced the notion of \emph{projection} on residues that  has been used extensively to study the abstract structure of buildings (see for instance \cite{RonanBuildings} or \cite{AbramBrown}).  Indeed, residues are sufficiently nicely embedded in $\Delta$ so that we can project the entire building on them  \emph{i.e} for any chamber $x \in \Delta$ and any residue $R\subset \Delta$ there exists a unique chamber $\proj{R}{x}\in R$ realizing the distance between $x$ and $R$.
 
 Two residues $R$ and $Q$ are   \emph{parallel} if  
 \[\proj{R}{Q} = \proj{Q}{R}.\] This notion has been introduced by J. Tits in \cite{TitsTwinBuildingsGroupsKacMoody} and is the object of an extensive study in \cite[Chapter 21]{MuhlherrPeterWeissDescentinBuildings}. These residues derive from  opposite residues in spherical buildings, with which they share a lot of properties. 
 
 The goal of this article is to study parallel residues and to relate this notion to residues with    equal stabilizers under a group action.

\subsection{Main results}
 
 With a simple geometric argument, we can observe that in a thin building parallelism is a transitive relation and thus is an equivalence relation on the set of residues. In the thick case, this holds if and only if $\Delta$ is right-angled (see   \cite[Proposition 2.10]{CapraceAutomRightAngled}). We will study the intermediate case and  characterize the buildings in which parallelism is a transitive relation by the structure of their  residues of rank 2. 
  
\begin{theo} [Theorem \ref{theoRelDequivCarac}]\label{theoprincip intro}
In a building $\Delta$, parallelism is an equivalence relation on the set of residues if and only if any  spherical residue of rank 2 is either thin or right-angled.
\end{theo}

 The second result is a group theoretical application of Theorem \ref{theoprincip intro}. Therefore we will
consider group actions on buildings. In this paper, all actions are assumed to be type-preserving. We recall that under this assumption, residues with equal stabilizers are parallel (see \cite[Proposition 22.3]{MuhlherrPeterWeissDescentinBuildings}). On the other hand, if the converse is true, then parallelism is an equivalence relation on the residues. In the thin and right-angled cases, the group is a Coxeter group or a graph product and we obtain the following corollary.

\begin{coro} [Corollary \ref{coroInterparagp}] 
Let $G$ be a Coxeter group or a graph product. Then, in $G$ intersections of parabolic subgroups are parabolic.
\end{coro}
  
   In the case of Coxeter groups, this Corollary is a classical fact due to J. Tits (see for instance \cite[Lemma 5.3.6]{DavisBook} for another proof). In the case of graph products this corollary has been established recently by Y. Antolín and A. Minasyan   by the means of Bass-Serre theory (see \cite[Proposition 3.4]{AntolinMinasyanTitsAlterGP}). The present article highlights in particular that these properties of Coxeter groups and graph products are true for the same reasons.

\subsection{Organization of the article}
In Section \ref{secBuildings}, we recall generalities about buildings, insisting on the notions of projections and  right-angled buildings. Then, in Section \ref{secParaRsidues} we discuss the notion of parallel residues and describe the buildings in which parallelism is an equivalence relation. Eventually, in Section \ref{secParaResiandStabi}, we study the situation where parallel residues admit the same stabilizers under a chamber-transitive group action.

\subsection{Terminology and notation}
\label{subsecterminota}

All along this article, we will use the following conventions. The identity element in a group will always be designated by $e$. For a set $E$,  the \emph{cardinality} of $E$ is designated by  $\# E$.  If $\mathcal{G}$ is a graph then $\mathcal{G}^{(0)}$ is the \emph{set of vertices} of $\mathcal{G}$ and $\mathcal{G}^{(1)}$ is the \emph{set of edges} of $\mathcal{G}$. For $v,w \in \mathcal{G}^{(0)}$, we write $v\sim w$ if there exists an edge in $\mathcal{G}$ whose extremities are $v$ and $w$.  
 
\subsection*{Acknowledgment} I am most grateful to the reviewer of this article for his careful reading. His fruitful comments  have been the origin of great improvements of this text. I would like to thank Pierre-Emmanuel Caprace for his advice and interest for this paper.  My thanks also go to Fanny Kassel for her careful reading of the first version of this article. Eventually, I thank my Ph.D advisor Marc Bourdon, for the  valuable and numerous discussions we had about this work.

\section{Buildings}
 \label{secBuildings}
Buildings are both combinatorial and geometric objects introduced by J. Tits to study Lie groups of exceptional types. In this section, we give a quick introduction to buildings. We emphasize the notion of a projection  and the particular case of right-angled buildings.  Buildings are the objects of extensive introductions in  \cite{RonanBuildings} and \cite{AbramBrown} to which we refer for details.
   \subsection{Chamber systems}

Throughout this paper $S$ is a fixed set.  

\begin{déf} A \emph{chamber system} $X$   over  $S$ is a set endowed with a family of partitions indexed by $S$. The elements of $X$ are called  \emph{chambers}.\end{déf} 
 
 In this subsection,  $X$ is a chamber system over $S$.  For $s\in S$, two chambers $c, c'\in X$ are said to be \emph{$s$-adjacent} if they belong to the same subset of $X$ in the partition associated with $s$. In this case, we write $c \sim_s c'$ and $s$ is called the \emph{type} of the adjacency relation. Usually, omitting the type  we refer to \emph{adjacent} chambers and we write $c \sim c'$. Note that any chamber is adjacent to itself.
 
 A map $f : X \longrightarrow X'$  between two chamber systems $X,X'$ over $S$ is a called a \emph{morphism} if it preserves the adjacency relations. If a morphism $f : X \longrightarrow X$ is a bijection, it is  called an \emph{automorphism} and  if moreover $f$ preserves the types of the adjacency relations, we say that $f$ is a   \emph{type preserving automorphism}. We designate by $\Autt{X}$ the group of  \emph{type preserving automorphisms} of $X$. Given a subset of $Y$ of $X$, then $Y$ inherits naturally the structure of a chamber system.

  We call   \emph{gallery}, a finite sequence  $\{c_k \}_{k=1, \dots , \ell}$  of chambers  such that  $c_k \sim c_{k+1}$  for  $k=1, \dots , \ell-1$. The galleries induce a metric on $X$.
   
  \begin{déf} \label{def distcombi} The \emph{distance between two chambers} $x$ and $y$ is the length of the shortest gallery connecting $x$ to $y$  and is designated by  $\dc{x, y}$. A shortest gallery   between two chambers is called \emph{minimal}.  \end{déf}
  
  For $I \subset S$, a subset $C$ of  $X$ is said to be \emph{$I$-connected} if for any pair of chambers $c,c' \in C$ there exists a gallery   $c=c_1 \sim \dots \sim c_\ell=c'$ such that for any $k=1, \dots , \ell-1$, the chambers  $c_k$ and $c_{k+1}$ are $i_k$-adjacent for some  $i_k \in I$. 
  \begin{déf} \label{def residus}
   The $I$-connected components are called the \emph{$I$-residues} or the \emph{residues of type $I$}. The \emph{rank} of an $I$-residue is the cardinality of $I$. The residues of rank $1$ are called \emph{panels}.  
   \end{déf}  
 
   We observe that a $I$-residue of a chamber system has a natural structure of a chamber systeme over $I$.
   
 A subset $C$ of $X$ is called \emph{convex} if every minimal gallery whose extremities belong to $C$ is entirely contained in $C$. 
Convexity is stable by intersection and for $A\subset X$,  the \emph{convex hull} of  $A$ is the smallest convex subset  containing $A$. In particular, convex subsets of $X$ are subsystems and residues are convex. 

The following example is crucial because it will be used to equip Coxeter groups and graph products with structures of chamber systems (see Definition \ref{def CoxSyst} and Theorem \ref{theo defbuildinggp}).  

\begin{exmp} \label{ex systeme chambre} Let  $G$ be a group, $B$ a subgroup and $\{H_i\}_{i\in I}$ a family of subgroups of $G$ containing $B$. The set of left cosets of $H_i/B$ defines a partition of $G/B$. We denote by $C(G,B,\{H_i\}_{i\in I})$ this chamber system over  $I$. This chamber system comes with a natural action of $G$. The group $G$ is a group of type-preserving automorphisms of $C(G,B,\{H_i\}_{i\in I})$ and the action is  chamber-transitive.\end{exmp}
In this paper we shall primarily be concerned with the case where $B=\{e\}$.

\subsection{Coxeter systems} 
 \label{subsect Coxeter}

 A  \emph{Coxeter matrix over $S$} is a symmetric matrix $M=\{m_{r,s}\}_{r,s\in S}$ whose entries are elements of $\N \cup \{\infty\}$ such that  $m_{s,s}=1$ for any $s\in S$ and  $\{m_{r,s}\}\geq 2$ for any $r, s \in S$ distinct. Let $M$ be a Coxeter matrix. The \emph{Coxeter group}   of type $M$ is the group given by the following presentation\[W=\left\langle s \in S \vert (rs)^{m_{r,s}}=1 \text{ for any } r,s\in S  \right\rangle .\]We call \emph{special subgroup} a subgroup of $W$  of the form \[W_I= \left\langle s \in I \vert (rs)^{m_{r,s}}=1 \text{ for any } r,s\in I  \right\rangle \text{ with } I\subset S. \]
 \begin{déf} \label{defparacox} A \emph{parabolic subgroup}    of $W$ is a subgroup  of the form $w W_I w^{-1}$ where $w\in W$ and $I\subset S$. An involution of the form $w s w^{-1}$ for $w\in W$ and $s\in S$ is called a \emph{reflection}.  \end{déf}

\begin{exmp} \label{ex cox poly}   Let $\mathbb{X}^d= \s^d, \mathbb{E}^d$ or $\h^d$. A \emph{Coxeter polytope}  is a convex polytope of $\mathbb{X}^d$ such that any dihedral angle  is of the form $\frac{\pi}{k}$ with $k$ not necessarily constant. Let $D$ be a Coxeter polytope and  let  $\sigma_1, \dots, \sigma_n$ be the codimension 1 faces of $D$.  We set $M= \{m_{i,j}\}_{i,j=1,\dots, n}$ the matrix defined by $m_{i,i}=1$, if $\sigma_i$ and $\sigma_j$ do not meet in a codimension 2 face  $m_{i,j}= \infty$, and if $\sigma_i$ and $\sigma_j$  meet in a codimension 2 face $\frac{\pi}{m_{i,j}}$ is the dihedral angle between $\sigma_i$ and $\sigma_j$. 

A theorem of  Poincaré (see for instance \cite[Theorem 1.2.]{GabPauImmeubles})  says that the reflection group  of $\mathbb{X}^d$  generated by the  codimension 1 faces of $D$ is  a discrete subgroup of $\mathrm{Isom}(\mathbb{X}^d)$ and is isomorphic to the Coxeter group of type $M$.
\end{exmp} 

\begin{déf} \label{def CoxSyst} With the notation introduced in Example \ref{ex systeme chambre}, the  \emph{Coxeter system} associated with $W$ is the chamber system over $S$ given by $C(W,\{e\},\{W_{\{s\}}\}_{s\in S})$. We use the notation $\coxsyst$   to designate this chamber system. \end{déf}

The chambers of $\coxsyst$ are the elements of $W$ and two distinct chambers $w,w' \in W$ are $s$-adjacent if and only if  $w=w's$. For $I\subset S$,  notice that  the  $I$-residues of $\coxsyst$ are the left-cosets of $W_I$ in $W$.  
 Again  $W$ is a group of automorphisms of $\coxsyst$ and the action is chamber-transitive.
 
Now we recall classical terminology about Coxeter systems.

\begin{déf}\numeroti{
\item Let $r=w s w^{-1}$  be a reflection for some $w\in W$ and $s\in S$. The  \emph{wall} $M_r$ in $\coxsyst$ is the  set of all the panels stabilized by $r$.
\item Let $M$ be a wall and $R$ be a residue. We say that $M$ \emph{crosses} $R$ if one of the panels of $M$ is contained in $R$. }
\end{déf}

In the particular case where $W$ is a finite group we refer to 	a \emph{spherical} Coxeter group and system. If $M=\{m_{r,s}\}_{r,s\in S}$ with $\{m_{r,s}\} \in \{2,\infty\}$ for any $r\neq s$, then we refer to a \emph{right-angled} Coxeter group or system. 
  
\subsection{Buildings} 
\label{subsec def build}

Hereafter $\cox$ is a fixed Coxeter system.
 \begin{déf} [{\cite[Definition 3.1.]{TitsBuildingsLectureNotes}}] A chamber system $\Delta$ over $S$ is a \emph{building of type $(W, S)$}   if it admits a maximal family $\apD$ of subsystems isomorphic to $\coxsyst$, called apartments, such that 

\begin{itemize}\item any two chambers lie in a common apartment,\item for any pair of apartments $A$ and $B$, there exists an isomorphism from $A$ to $B$  fixing $A\cap B$.\end{itemize}\end{déf}

   If the group $W$ is a spherical (resp. right-angled) Coxeter group then $\Delta$ is called a \emph{spherical} (resp. \emph{right-angled}) building. 

 Hereafter, $\Delta$ is a fixed building of type $\cox$. A straightforward application of this definition is the existence of retraction maps of the building over apartments.\begin{déf}Let $x\in \Delta$ and $A\in \apD$. Assume that $x$ is contained in $A$. We call   \emph{retraction onto $A$ centered $x$} the map $\pi_{A,x} : \Delta \longrightarrow A$ defined by the following property. \begin{center}\begin{minipage}[c]{13cm}For $c\in \Delta$, there exists a chamber $\pi_{A,x}(c) \in A$ such that for any apartment $A'$  containing $x$ and $c$,  for any  isomorphism  $f:A'\longrightarrow A$  that fixes $A \cap A'$, then $f(c)=\pi_{A,x}(c)$\end{minipage}\end{center}\end{déf}

\begin{exmp} \label{ex building}
 \numeroti{
 \item Any infinite tree without leaf is a building of type $(W,S)$ where $W$ is the infinite dihedral group $\Z/2\Z*\Z/2\Z$ and $S= \{(1,0), (0,1)\}$.
 
 \item 
For $n\geq 1$ and $\mathbf{k}$ a field, the flags of  subspaces of a $n$ dimensional vector space over $\mathbf{k}$ is  a spherical building (see \cite[Chapter 1]{RonanBuildings}). On   Figure \ref{fig:D0} is represented the geometric realisation of the building of $\mathbf{k}^{3}$ where  $\mathbf{k}$ is the finite field of order $2$. 
}   
\end{exmp} 
 
 The building $\Delta$ is called a \emph{thin} (resp. \emph{thick})   building if any panel contains exactly two (resp. at least three) chambers. Note that thin buildings are Coxeter systems.
 We recall that an $I$-residue of $\Delta$ is itself a building of type $(W_I,I)$. Hence it makes sense to talk about \emph{thin}, \emph{thick}, \emph{spherical} or \emph{right-angled} residues.

For  $x$ and $y$ two chambers of $\Delta$, the  convex hull of the pair $\{x,y\}$ in $\Delta$ is the convex hull of   $\{x,y\}$ in any apartment containing $x$ and $y$ (see  \cite[Proposition 3.18.]{TitsBuildingsLectureNotes}). This fact permits to build projections on residues. 

\begin{prop}[{\cite[Proposition 3.19.3.]{TitsBuildingsLectureNotes}}] \label{prop defprojection} Let $R$ be a residue and $x$ be a chamber in $\Delta$.  There exists a unique chamber  $\proj{R}{x}\in R$ such that $\dc{x,\proj{R}{x}}= \di{x, R}$.   Moreover, for any chamber  $y$ in  $R$ there exists a minimal gallery from  $x$ to $y$ passing through   $\proj{R}{x}$. \end{prop}
 
 Observe that not all convex subsets of a building admit projection maps. Indeed, let $P$ be a panel of $\Delta$. Then any subset $C$ of $P$ is convex. However, the projection of $\Delta$ onto $C$ exists if and only if $\#C = 1$ or $C = P$.
 
\subsection{Graph products and right-angled buildings}

Let $\mathcal{G}$ denote a \emph{simplicial graph} \emph{i.e} no edge  is a loop and no edge is double.  If in $\mathcal{G}$ two distinct  vertices $v$ and $v'$ are connected by an edge, we write $v \sim v'$.  A   group $ G_v$    is associated with each $v \in \mathcal{G}^{(0)}$ and we denote by $F_\mathcal{G}$ the free product of the family $\{G_v\}_{v\in \mathcal{G}^{(0)}}$.
 
\begin{déf} \label{def graphprod}The \emph{graph product}   given by the pair $(\mathcal{G},\{G_v\}_{v\in \mathcal{G}^{(0)}})$ is the group defined by the following quotient \[\Gamma = F_\mathcal{G}/ R, \]
where $R$ is the normal subgroup $\left\langle\left\langle  gg'g^{-1}g'^{-1} : g\in G_v, g'\in G_{v'} \text{ and } v\sim v' \right\rangle\right\rangle$. \end{déf}

\begin{exmp} If all the groups $\{G_v\}_{v\in \mathcal{G}^{(0)}}$ are of order $2$ then $\Gamma$ is  a right-angled Coxeter group. If all the groups $\{G_v\}_{v\in \mathcal{G}^{(0)}}$ are infinite cycles, then  $\Gamma$ is a   \emph{right-angled Artin group} (see \cite{CharneyRAAG}). In fact, all right-angled Coxeter and Artin groups may be obtained as a graph product. If $\mathcal{G}$ has no edge $\Gamma$ is a free product and if  $\mathcal{G}$  is a complete graph $\Gamma$ is a free Abelian product.\end{exmp}
 
Now  we designate by $S$ the set $\mathcal{G}^{(0)}$. This is motivated by the fact that a Coxeter group is canonically associated to a graph product. From now on, we fix a graph product $\Gamma$ given by a pair  $(\mathcal{G},\{G_s\}_{s\in S})$. Then,  the graph product  defined by the pair $(\mathcal{G},\{\Z/2\Z \}_{s\in S})$ is isomorphic to the right-angled Coxeter group defined by the matrix $M=\{m_{s,t}\}_{s,t \in S}$ given by: $m_{s,t}=2$ if $s\sim t$ and  $m_{s,t}=\infty$ if $s\nsim t$ in $\mathcal{G}$. We denote by $W$  this Coxeter group and by $\cox$ the associated Coxeter system.

With this notation, the following theorem associates a right-angled building to  a graph product.  
 
\begin{theo}[{\cite[Theorem 5.1.]{DavisCAT}}] \label{theo defbuildinggp}
Let $\Delta$ be the chamber system $C(\Gamma,\{e\}, \{G_{s}\}_{s\in S})$  (see Example  \ref{ex systeme chambre}). Then $\Delta$ is a building of type $\cox$.
\end{theo}
 
 A classification of F. Haglund and F. Paulin  states that the construction presented above describes all the right-angled buildings in which all the panels of same type have same cardinality.
 
\begin{theo}[{\cite[Proposition 5.1.]{HaglundPaulinImmeubles}}]\label{theoClassifRAB} Let $\Gamma$ be the graph product given by the pair $(\mathcal{G},\{G_s\}_{s\in S})$. Let $\Delta$ be the building of type $\cox$ associated with  $\Gamma$ by Theorem \ref{theo defbuildinggp}.  Assume that  $\Delta'$ is a building of type $\cox$ such that for any $s\in S$ the  $\{s\}$-residues of $\Delta'$ are of cardinality $\#G_s$. Then $\Delta$ and $\Delta'$ are isomorphic. 
\end{theo}  
 
 By analogy with Definition \ref{defparacox},  we define parabolic subgroups in $\Gamma$.
 \begin{déf} \label{defparagp}
For $I\subset S$ we write $\Gamma_I=\left\langle G_s : s\in I\right\rangle$ and   a subgroup of the form $g \Gamma_I g^{-1}$, with $g \in \Gamma$, is called a \emph{parabolic subgroup} of $\Gamma$.  
    \end{déf}
      
\section{Parallel residues} 
  \label{secParaRsidues}
  
Parallel residues have been defined by J. Tits in \cite{TitsTwinBuildingsGroupsKacMoody}. This notion derives from the notion of opposite residues.  We refer to \cite[Chapter 21]{MuhlherrPeterWeissDescentinBuildings} for details about parallel residues in general and to \cite[Chapters 5 and 9]{WeissSphericalBuildings} for details about opposite residues. 

The goal of this section is to study the buildings in which parallelism of residues is a transitive relation.

\subsection{Definition and first properties}

In the rest of the paper, $\Delta$ is a building of type $\cox$.  

\begin{déf} \label{def paralresidu}
Let $R$ and $Q$ be two residues in $\Delta$. We say that $R$ is \emph{parallel} to $Q$ if 
\[ \proj{R}{Q}=R \text{ and } \proj{Q}{R}=Q.\]
\end{déf}

The following proposition summarizes some basic properties of parallel residues.

\begin{prop}[{\cite[Propositions 21.8 and 21.17]{MuhlherrPeterWeissDescentinBuildings}}]\label{propBaseResiduesPara}
Let $R$ and $Q$ be respectively a  $I$-residue  and a $J$-residue in $\Delta$. Let  $Q'=\proj{Q}{R}$ and $R'=\proj{R}{Q}$. Then the following properties hold.
\numeroti{
\item $R$ is parallel to $Q$ if and only if for any apartment $A$ containing a chambers of both $R$ and $Q$ the residues $R\cap A$ and $Q\cap A$ are parallel in $A$. 
\item $R'$ and $Q'$ are parallel residues.
\item   The maps ${\mathrm{proj}_R}_{\vert Q'} : Q'\longrightarrow  R'$ and  ${\mathrm{proj}_Q}_{\vert R'} :R'\longrightarrow Q'$ are  reciprocal  bijections. 
\item For any $x,y \in R'$,  $\dc{x,\proj{Q}{x}}= \dc{y,\proj{Q}{y}}$

\item There exists a unique $w{(R,Q)} \in W$ such that for any apartment $A$ containing a chambers of both $R'$ and $Q'$, for any chamber  $x$ in $R'\cap A$ one has in $A$ : $w{(R,Q)}x=\proj{Q'}{x}$.

\item Let $w = w{(R,Q)}$, then  $R'$ (resp. $Q'$) is of type $I'= \{s\in I : w^{-1} s w = t \text{ for some } t \in J \}$ (resp. $J'= \{s\in J : w^{-1} s w = t \text{ for some } t \in I \}$).}

\end{prop}

\subsection{Opposite residues}

In this subsection,  $\Delta$ is a spherical building. If not further specified, the proofs of the following claims are contained in  \cite[Chapters 5 and 9]{WeissSphericalBuildings}. 

For a chamber $x$ in  $\Delta$, a chamber $y$ is called \emph{opposite} to $x$ if $\dc{x,y}=\dia{\Delta}$, where $d_c$ is the distance over the chambers given by Definition \ref{def distcombi}. Clearly this definition is empty in the non-spherical case. However, it is very rich in the spherical case. Indeed,   for any chamber $x$ and apartment $A$ containing $x$, there exists a unique   chamber  $y$ opposite to $x$ contained in $A$. We denote by $\mathrm{op}_A : A \longrightarrow A$ the map sending a chamber to its opposite chamber in $A$. 

\begin{déf}
Let $R$ and $Q$ be two residues in $\Delta$. We say that $R$ and  $Q$ are \emph{opposite} residues if there exists an apartment $A$ intersecting both $R$ and $Q$ such that \[ \mathrm{op}_A(R) = Q\cap A \text{ and }  \mathrm{op}_A(Q)=R\cap A.\]
\end{déf}

In fact two residues are opposite if and only if the condition of the preceding definition is satisfied for any apartment $A$ intersecting both $R$ and $Q$.

Opposite residues are parallel (see for instance\cite[Proposition 21.24]{MuhlherrPeterWeissDescentinBuildings}).  However the converse is false as any two chambers are always parallel. In fact, the notions of opposite and parallel residues in a spherical building are connected by the following proposition.

\begin{prop}[{\cite[Proposition 21.26]{MuhlherrPeterWeissDescentinBuildings}}]
Two parallel residues $R$ and $Q$ are opposite if and only if for some chamber $x\in R$ there exists a chamber $y\in Q$ opposite to $x$ in $\Delta$. 
\end{prop}

Moreover, it appears that for any residue $R\subsetneq \Delta$ and apartment $A$, there exists a residue $Q$ such that $A$ intersects $Q$ and $R$ is opposite to $Q$.  Thus, any residue admits an opposite (and thus a parallel) residue. This is not true in the non-spherical case.  Indeed, in the thin building associated with the infinite dihedral group, no panel admits a parallel residue.

\subsection{Parallel residues defining an equivalence relation}

By simple geometric arguments, we observe that in a thin building parallel residues are characterized by walls.

\begin{prop}[{\cite[Proposition 21.19]{MuhlherrPeterWeissDescentinBuildings}}] \label{Propmurcoxeter} In a thin building, two residues are parallel if and only if the set of walls that cross them are equal.\end{prop}

\begin{figure}[h!]
    \labellist
\small\hair 2pt
\pinlabel $T$ at 35 425
 
\pinlabel $Q$ at 230 570 
 
\pinlabel $R$ at 355 -10

\endlabellist
\centering
\includegraphics[scale=0.254]{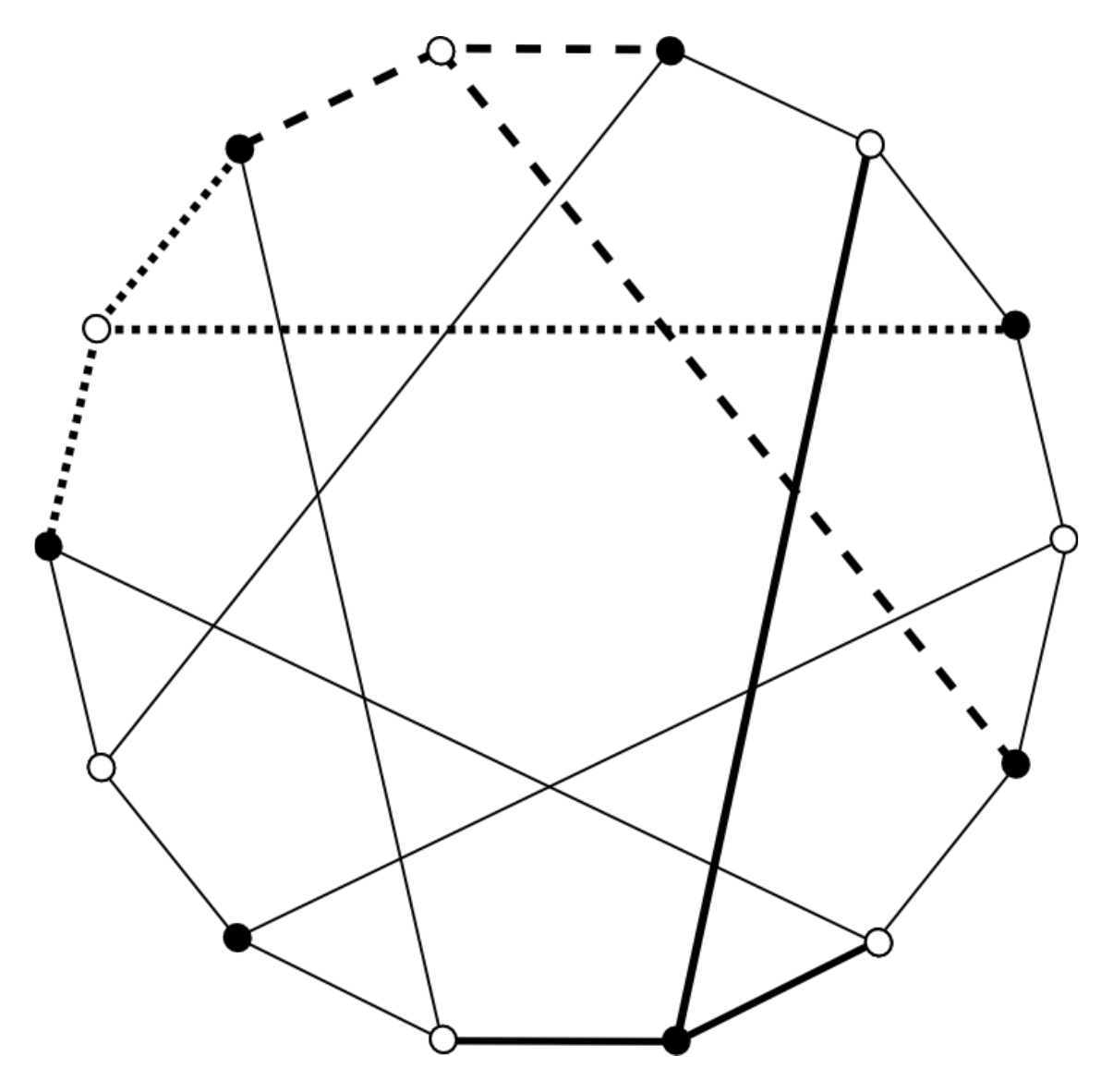}
\caption{In this spherical building $R$ is parallel to $Q$ and $T$ but $Q$ is not parallel to $T$.}
\label{fig:D0}
\end{figure}

 A consequence of the preceding proposition is that  in a thin building, parallelism is an equivalence relation  on  the residues.  
In general, the relation induced by the  parallelism may not be  transitive as illustrated by the example of Figure \ref{fig:D0}.  
 In fact, in the thick case,  this happens if and only if the building is right-angled (see   \cite[Proposition 2.10]{CapraceAutomRightAngled}).  

  In the following we observe that this strong property leaves only few examples between thin and right-angled buildings. To this end we will use several times the following  fact: two panels $\sigma$ and $\sigma'$ are parallel if and only if there exists an apartment $A$ in which $\sigma\cap A$ is parallel to $\sigma'\cap A$. This follows directly from the definition of the projections.

We start by establishing a short lemma about thin and right-angled residues in $\Delta$.

\begin{lem} \label{lem PaneauFinRA}
Let $R$ be a thin or a right-angled residue in $\Delta$ , let $\sigma$ and $\sigma'$ be two parallel panels in $R$ and let  $\delta$ be a panel in $\Delta$.  If   $\delta$ is parallel to $\sigma$ then it is parallel to $\sigma'$. 
\end{lem}

\begin{proof}
If $R$ is thin, then any apartment containing $\sigma$ also contains $\sigma'$. Hence the lemma is satisfied by Propositions \ref{propBaseResiduesPara}.$i)$ and  \ref{Propmurcoxeter}. 
 
Now we assume that  $R$ is right-angled. We define  $\delta':= \proj{R}{\delta}$ and we  observe that as   $\sigma$ is contained in $R$ one has $\proj{\sigma}{\delta}=\proj{\sigma}{\delta'}$. If $\delta'$ is not parallel to $\sigma$ then  $\proj{\sigma}{\delta'}$ is a single chamber. But this is absurd because  $\proj{\sigma}{\delta'} = \proj{\sigma}{\delta}$ and $\sigma$ is parallel to $\delta$. Likewise, if $\delta'$ is not parallel to $\delta$ then $\delta'$ is a single chamber. This implies again that  $\proj{\sigma}{\delta}$ is a single chamber which is absurd. Hence $\delta'$ is parallel to both $\sigma$ and $\delta$.

 As in $R$  parallelism is an equivalence relation, $\delta'$ is parallel to $\sigma'$. By convexity of the apartments and by Proposition \ref{prop defprojection}, any apartment intersecting both $\delta$ and $\sigma'$ intersects $\delta'$. Then, by transitivity of parallelism in the apartments  and by  Propositions \ref{propBaseResiduesPara}.$i)$, we obtain that $\delta$ is parallel to $\sigma$.
\end{proof}

To characterize buildings in which parallelism is an equivalence relation on the set of residues, we will use the following notation. For $\sigma$ and $\delta$  two panels in $\Delta$, for an apartment $A$ intersecting both $\sigma$ and $\delta$ we write $r$ and $t$ for the reflections in $A$ stabilizing respectively $\sigma\cap A$ and $\delta\cap A$. We call \emph{order} of the pair $(\sigma,\delta)$, and we write $\mathrm{Ord}(\sigma,\delta)$, the order of $rt$ in $W$ and we observe that   $\mathrm{Ord}(\sigma,\delta)$ is well defined \emph{i.e} it does not depend on the choice of $A$. Moreover, as a consequence of Proposition \ref{propBaseResiduesPara}.$i)$, $\sigma$ is parallel to $\delta$, if and only if $\mathrm{Ord}(\sigma,\delta)=1$.
   
\begin{theo} \label{theoRelDequivCarac}  In a building $\Delta$, the following are equivalent:
\numeroti{\item  Parallelism  is an equivalence relation on the set of  residues.
 \item Parallelism is an equivalence relation on the set    of the panels.
\item  For any pair of panels $\sigma$ and $\delta$, if $\sigma$ is thick then $\mathrm{Ord}(\sigma,\delta)\in \{1,2, \infty\}$.
\item Any  spherical residue of rank $2$ is either thin or right-angled.
}
\end{theo}

\begin{proof}
The implications $i) \Longrightarrow ii)$ and $iii) \Longrightarrow iv)$ are immediate and we start by proving $ii) \Longrightarrow i)$. Let $R$, $Q$ and $T$ be residues such that $R$ and $Q$ are parallel to $T$. Let $A$ be an apartment intersecting both $R$ and $Q$. Let $M$ be a wall crossing $A\cap R$. By Propositions \ref{propBaseResiduesPara}.$i)$ and \ref{Propmurcoxeter}, it is enough to prove that $M$ crosses $A\cap Q$. To this end let $\sigma\subset R$ be a panel such that   $\sigma\cap A$ is  crossed by $M$.  Then, as $R$ is parallel to $T$, there exists a panel $\sigma_T$ in $T$ that is parallel to  $\sigma$. 

Now pick an apartment $A'$ intersecting both $\sigma_T$ and $A\cap Q$.  As $A'\cap T$ is parallel to $A'\cap Q$, there exists a panel $\delta_{A'}$ in $A' \cap Q$ that is parallel to  $\sigma_T\cap A'$. We observe that the panel $\delta\subset Q$ containing $\delta_{A'}$ is parallel to $\sigma_T$. Then, by  assumption, $\sigma$ is parallel to $\delta$ and by Proposition \ref{Propmurcoxeter}, $M$ crosses $\delta\cap A$. 
 
We prove  $ii) \Longrightarrow iii)$ by contradiction (this step is essentially the same as the proof of \cite[Proposition 2.10]{CapraceAutomRightAngled}). Let  $\sigma$ be a thick panel and $\delta$ be a panel such that $\mathrm{Ord}(\sigma,\delta)=n>2$. In an apartment $A$ intersecting  both $\sigma$ and $\delta$, the wall crossing $\sigma\cap A$ intersects the wall crossing $\delta\cap A$. As a consequence, $\Delta$ contains a residue $R$ of rank 2 that is not right-angled nor thin. 

 Then, we set $\sigma'=\proj{R}{\sigma}$ and we choose two distinct panels $\sigma_1$, $\sigma_2$ of the same type, contained in $R$, lying at a minimal distance and containing a chamber of $\sigma'$.  Choose an apartment $A$ and a chamber $x\in A$ such that $\pi_{A,x}(\sigma_1)= \pi_{A,x}(\sigma_2)$. In $R\cap A$, let $\delta'_A$ be the panel opposite to $\pi_{A,x}(\sigma_1)$ and let  $\delta'$ be the panel in $\Delta$ such that $T\cap A = T_A$. As $\pi_{A,x}$ decrease the distance over the chambers, then  $\delta'$ is opposite to both $\sigma_1$ and $\sigma_2$ in $R$ and thus is parallel to them in $\Delta$.
 
Here we prove $iv) \Longrightarrow ii)$. Let $\sigma$, $\sigma'$ and $\delta$ be three panels such that $\sigma$ and $\sigma'$ are parallel to  $\delta$. We prove the implication by induction on $d=\max\{ \di{\sigma,\delta} ,\di{\sigma',\delta}\}$.  If $d=0$ there is nothing to prove. 

If $d>0$, consider $R$ a  residue of rank $2$ containing $\sigma$ and such that  $\di{\delta,R}< \di{\delta,\sigma}$ and choose an apartment $A$  intersecting both $\sigma$ and $\delta$. There exists a panel $T_A$ in $R\cap A$ such that $T_A$ is parallel to $\sigma\cap A$ and   $\di{\delta,R} =  \di{\delta\cap A, T_A}$. This panel is the panel opposite to $\sigma\cap A$ in  $R\cap A$. We designate by $T$ the panel in $R$ containing $T_A$  and we check that $T$ is parallel to both $\sigma$ and $\delta$. We do the same with   $\sigma'$ and we obtain $T'$  parallel to both $\sigma'$ and  $\delta$ and such that  $\di{T',\delta}<d$. 

Now, by the induction assumption, we obtain that $T$ is parallel to $T'$. To finish, we observe that $\sigma$ and $T$ (resp. $\sigma'$ and $T'$)  are contained in thin or right-angled residues and the proof is achieved  by  Lemma \ref{lem PaneauFinRA}.   \end{proof}

As it is suggested by the preceding theorem, the buildings in which parallelism is an equivalence relation are obtained from right-angled buildings by substituting a given Coxeter system for chambers. Here we explain this fact in detail.

In the rest of this section,  $\Delta$ is a building of type $\cox$  satisfying  the equivalent conditions of Theorem \ref{theoRelDequivCarac}. Let $M= \{m_{s,r}\}_{s,r \in S}$ be the Coxeter matrix associated to $\cox$. We set:  \liste{\item $S_\perp:=\{s\in S: m_{s,r}\in\{2, \infty \} \text{ for any } r\neq s\}$, 
\item $S_T:=S\backslash S_\perp $.} The set $S_\perp$ is the set of possibly thick  types of $\Delta$. We designate by  $R_T$ a  thin residue in $\coxsyst$ of type $S_T$ and we define the following graph $\mathcal{G}$: 

\begin{itemize}
\item $\mathcal{G}^{(0)} = \{wsw^{-1}\in W : w\in W_{S_T} \text{ and } s\in S_\perp \}$. Equivalently,  $\mathcal{G}^{(0)}$ is the set of walls that bound $R_T$ in $\coxsyst$. 
\item Two vertices $v,v' \in \mathcal{G}^{(0)}$ are joined by an edge if and only if  the corresponding reflections commute.  Equivalently, if and only if  the corresponding walls intersect in $\coxsyst$. 
\end{itemize}

Now we designate by $S'_\perp$ the set  of vertices of $\mathcal{G}$  and by $(W_\perp, S'_\perp)$ the Coxeter system associated to $\mathcal{G}$. By construction, it appears that the set of all the $S_T$-residues  of $\Delta$  inherits from   $\Delta$   a structure of right-angled building of type $(W_\perp, S'_\perp)$. 
We denote by $\Delta_\perp$ this building. Observe that $S_\perp$ is not always equal to $S'_\perp$. For instance if $S$ is finite and if there exists $s\in S_\perp$ and $r,t\in S_T$ such that $m_{s,t}=m_{r,t} = \infty$, then $S'_\perp$ is infinite. 

From now on, we assume that in   $\Delta_\perp$ all panels of the same type are of the same cardinality. For each $s\in S'_\perp$ we fix a group $G_s$ such that  $\#G_s= \#\sigma_s$.
  Then, by Theorem \ref{theo defbuildinggp}, $\Delta_\perp$ is isomorphic to the right-angled building associated to the graph product $\Gamma$ given by the pair $(\mathcal{G}, \{G_s\}_{s \in S'_\perp})$.   In particular, $\Gamma$ acts on  $\Delta$  with quotient  equal to $R_T$ and in fact, under these assumptions, $\cox$ and $\Gamma$ determine $\Delta$ up to isomorphism.

\begin{prop} Let $\Delta$ be a building of type $\cox$ satisfying the equivalent conditions of Theorem \ref{theoRelDequivCarac} and let $\Delta_\perp$ be the right-angled building of the $S_T$-residues of $\Delta$. If in   $\Delta_\perp$ all panels of the same type are of the same cardinality, then $\Delta$  is uniquely determined, up to isomorphism, by these cardinalities. 
\end{prop}

\begin{proof}
 Let $\Delta$ and $\Delta'$ be two buildings of same type $\cox$ satisfying  the equivalent conditions of Theorem \ref{theoRelDequivCarac}. Let $\Delta_\perp$ and $\Delta'_\perp$ be the two right-angled buildings  associated to them and assume that in $\Delta_\perp$ and $\Delta'_\perp$  all panels of the same type are of the same cardinality.  If these cardinalities are equal then $\Delta_\perp$ and $\Delta'_\perp$ are both isomorphic to the right-angled building given by the graph product $\Gamma$ defined as in the preceding paragraph.

Now we fix two base chambers $x_0\in \Delta$ and $x'_0\in \Delta'$. We designate by $R_T$ and $R'_T$ the  $S_T$-residues containing  respectively $x_0$ and $x'_0$ and we consider the isomorphism $f:R_T \longrightarrow R'_T$ mapping $x_0 \longmapsto x'_0$.  Then we observe that $f$ extends as a building isomorphism $F:\Delta \longrightarrow\Delta'$ as follow. For $x\in \Delta$, let $g\in \Gamma$ be such that $x=\gamma y$ with $y\in R_T$ then \[F(x): = \gamma f(y).\]
\end{proof} 

As a particular case we obtain the following corollary.

\begin{coro}
 Let $\Delta$ be a building of type $\cox$ satisfying the equivalent conditions of Theorem \ref{theoRelDequivCarac}. If in  $\Delta$ all panels of the same type are of the same cardinality, then $\Delta$  is uniquely determined, up to isomorphism, by these cardinalities. 
\end{coro}
 
 \section{Parallel residues and stabilizers}
 \label{secParaResiandStabi}
In this section, $G$ is a subgroup of $\Autt{\Delta}$ acting chamber-transitively. Here we discuss the relationship between the fact that two residues are parallel and the fact that these two residues have same stabilizers under the action of $G$.

 \subsection{Parallel residues with equal stabilizers}
  First, we  recall that as the action of $G$ is chamber-transitive, then    two residues with equal  stabilizers  are parallel (see \cite[Proposition 22.3]{MuhlherrPeterWeissDescentinBuildings}).  The following proposition describes the situation where the converse is true.

 \begin{prop} \label{PropPrincip} Suppose that $\mathrm{Stab}_G(P)=\mathrm{Stab}_G(P')$ for any pair  $P,P'$ of parallel residues. Then the following properties are satisfied: 
 
 \numeroti{\item Parallelism is an equivalence relation on the residues. In particular, $\Delta$ satisfies the equivalent conditions of Theorem \ref{theoRelDequivCarac}.
 
 \item The action is free.

 \item For any pair of residues $R,Q$ one has 
\[ \mathrm{Stab}_G(R)\cap\mathrm{Stab}_G(Q) = \mathrm{Stab}_G(\proj{R}{Q}) =\mathrm{Stab}_G(\proj{Q}{R}) .\] }
 
\end{prop}

\begin{proof}

$i)$ Under the hypothesis of the proposition, two residues have same stabilizer if and only if they are parallel.

$ii)$ Let $x$ be a chamber in $\Delta$. As any pair of chambers are parallel residues, for all $y\in \Delta$  one has $\mathrm{Stab}_G(x)=\mathrm{Stab}_G(y)$. Thus $\mathrm{Stab}_G(x)=\{e\}$

$iii)$ By symmetry, it is sufficient to prove, that  \[\mathrm{Stab}_G(R)\cap\mathrm{Stab}_G(Q)=\mathrm{Stab}_G(\proj{Q}{R}).\]Let $g \in \mathrm{Stab}_G(R)\cap\mathrm{Stab}_G(Q)$. As $g$ is an automorphism of $\Delta$ that stabilizes both $R$ and $Q$, the map $\proj{Q\vert R}{\cdot}$ is equivariant by $g$. Then $g(\proj{Q}{R})=\proj{Q}{R}$ and \[\mathrm{Stab}_G(R)\cap\mathrm{Stab}_G(Q)<\mathrm{Stab}_G(\proj{Q}{R}).\]Let $g \in \mathrm{Stab}_G(\proj{Q}{R})$. As $g$ preserves the types, if $Q$ is a $I$-residue then $g(Q)$ is also a $I$-residue. In particular, $Q$ and $g(Q)$ are two $I$-residues containing $\proj{Q}{R}$, thus $g(Q)=Q$. As $\proj{Q}{R}$ is   parallel to $\proj{R}{Q}$, under our assumption $g\in  \mathrm{Stab}_G(\proj{R}{Q})$. We can use the previous argument to prove that $g(R)=R$ and \[\mathrm{Stab}_G(\proj{Q}{R})<  \mathrm{Stab}_G(R)\cap\mathrm{Stab}_G(Q).\]

\end{proof}
  
 In the rest of the section, we   assume that the action of $G$  is chamber-transitive and that the assumption of the preceding proposition hold. 
 
  In the thin case, it is clear that $G$ is isomorphic to $W$.   In the right-angled case, the next proposition says that it is isomorphic to a graph product of stabilizers of panels.  To this end, we will use the following notation. For a right-angled Coxeter group $W$, we designate by  $\mathcal{G}_W$ the simplicial graph such that the graph product  given by $(\mathcal{G}_W,\{\Z/2\Z\}_{s\in S})$ is isomorphic to $W$. 
 
 \begin{prop} \label{proptransitivDcGP}Let $\Delta$ be a right-angled building of type $\cox$ and $G$ a group of type preserving automorphisms acting freely and chamber-transitively on $\Delta$.  Let $x_0$ be a chamber in $\Delta$ and let $G_s$ be the stabilizer in $G$ of the $s$-panel containing $x_0$. Then $G$ is isomorphic to the graph product given by the pair $(\mathcal{G}_W,\{G_s\}_{s\in S})$.
  \end{prop}
  
\begin{proof}
Let  $\Gamma$ be the graph product given by the pair $(\mathcal{G}_W,\{G_s\}_{s\in S})$. To $\Gamma$ we associate the right-angled building $\Delta_\Gamma$ given by Theorem \ref{theo defbuildinggp}. We observe that $\Delta_\Gamma$ is of type $\cox$ and that for  $\sigma_s(\Delta)$ and $\sigma_s(\Delta_\Gamma)$ two panels of type $s\in S$ respectively in $\Delta$ and in $\Delta_\Gamma$ one has:
\[ \#\sigma_s(\Delta)= \#\sigma_s(\Delta_\Gamma).\]

Hence, by Theorem \ref{theoClassifRAB},  $\Delta$ and $\Delta_\Gamma$ are isomorphic and we both denote them $\Delta$. As a consequence, $\Gamma$ is the subgroup of $\Autt{\Delta}$   generated by the set $\{G_s\}_{s\in S}$. In particular, this proves that $\Gamma < G$.

Now, for $g\in G$, we prove  by induction on $n=\dc{x_0,gx_0}$  that $g$ is a product of elements of   $\{G_s\}_{s\in S}$. If $n=0$ there is nothing to prove. If $n>0$ consider a minimal gallery: \[x_0\sim  \dots \sim x_{n-1} \sim x_n=gx_0.\]

Let $h\in G$ be such that $hx_{n-1} =x_n$. As $h$ preserves the type, $h\in \mathrm{Stab}_G(\sigma)$ where $\sigma$ is the $s$-panel containing $\{x_{n-1}, x_n\}$. Let $\gamma \in G$ be such that $\gamma x_0=x_{n-1}$. In particular, $\sigma= \gamma \sigma_s$,    $\mathrm{Stab}_G(\sigma)=\gamma G_s \gamma^{-1}$ where $\sigma_s$ is the $s$-panel containing $x_0$ and $h= \gamma g_s \gamma^{-1}$ for one $g_s\in G_s$. Then, by freeness of the action, $g=h\gamma$   and  with $ \di{x_0,\gamma x_0}   = n-1$  the proof is achieved.
 
\end{proof}

\subsection{Application to intersection of parabolic subgroups}

In this section, we apply Proposition \ref{PropPrincip} to thin and right-angled buildings under the action of Coxeter groups and graph products.

First we verify that the assumption of the theorem are satisfied in the case of a Coxeter groups.
 
\begin{prop}
If $\Delta$ is a thin building, then parallel residues have equal stabilizers. \end{prop}

\begin{proof}

Let $R$ be a residue. Here we prove that the  stabilizer of $R$ under the action of $W$ is the subgroup $G<W$ generated by the reflections about the walls that cross $R$. This will imply the proposition by  Proposition \ref{Propmurcoxeter}. 

As $W$ is type preserving, it is clear that $G<\mathrm{Stab}_W(R)$. Now we fix $x_0\in R$ and for $g\in \mathrm{Stab}_W(R)$ we consider a minimal gallery
\[x_0\sim x_1\sim\dots\sim x_n=gx_0.\]
By convexity of the residues, this gallery is contained in $R$. Let $r_i\in W$ be the reflection that maps $x_{i}$ to $x_{i+1}$. Then, by simple chamber-transitivity of the action, $g= r_n \dots r_0$ and the proof is complete.
\end{proof}

In the right-angled case we establish an analogue proposition.

  \begin{prop}  Let $\Gamma$ be the graph-product given by a pair $(\mathcal{G},\{G_s\}_{s\in S})$ and let $\Delta$ be the associated right-angled building.  Then any two parallel residues of $\Delta$ have equal stabilizers.\end{prop}
  \begin{proof}
  Let $R$ and $Q$ be two parallel residues. Up to a conjugation, we can assume that $x_0$ is in $R$. According to Proposition \ref{propBaseResiduesPara}.$vii)$, $R$ and $Q$ are of same type $I$. We write
  \[ I^{\perp} = \{s \in S \backslash I : v_s \sim v_i \text{ for all } i\in I\}.\]
 By \cite[Proposition 2.8.$ii)$]{CapraceAutomRightAngled}, $R$ and $Q$ are both contained in $T$ a $J$-residue where $J=I\cup I^\perp$. We observe that $ \Gamma_J=\Gamma_I\times \Gamma_{I^\perp}$ and that $\mathrm{Stab}_\Gamma(T)=\Gamma_J$. As a consequence, $\Gamma_J$ acts transitively on the set of $I$ residues contained in $T$. Thus, there exists $g\in \Gamma_J$ such that $gR=Q$. Hence $\mathrm{Stab}_\Gamma(Q)=g \mathrm{Stab}_\Gamma(R) g^{-1}$ and with $\mathrm{Stab}_\Gamma (R)= \Gamma_I$ the proposition is proved.
  \end{proof}
  
    Now we know that both actions of Coxeter groups and of graph-products on their associated buildings satisfy the assumption of Proposition \ref{PropPrincip}. In the next proposition we obtain from this fact that intersections of parabolic subgroups are parabolic.

    From now on, $\Delta$ is either a thin or a right-angled building of type $\cox$. We fix a base chamber $x_0\in \Delta$ and  for $s\in S$ we denote by $\sigma_s$ the $s$-panel containing $x_0$. The group $G$ is a group acting freely and chamber-transitively on $\Delta$. In fact, $G$ is either $W$ in the thin case or a graph product $\Gamma$ in the right-angled case (see Proposition \ref{proptransitivDcGP}).  For $I\subset S$ we set \[G_I:= \left\langle \mathrm{Stab}_G(\sigma_s) :s\in I \right\rangle .\]
In fact,  $G_I$ is either $W_I$ in the thin case or $\Gamma_I$ in the right-angled case (see Definitions \ref{defparacox} and \ref{defparagp}). We recall that  a parabolic subgroup $gG_I g^{-1} < G$ stabilizes the $I$-residue $R=gG_I x_0$. We also recall that, according to Proposition \ref{propBaseResiduesPara}.$v)$,   for $R$ and $Q$ two residues, $w{(R,Q)} \in W$ is  such that for any apartment $A$ containing a chambers of both $\proj{R}{Q}$ and $\proj{Q}{R}$ and for any chamber  $x$ in $\proj{R}{Q}\cap A$ one has in $A$: $w{(R,Q)}x=\proj{Q}{x}$.       
 \begin{coro} \label{coroInterparagp} For $g\in G$, and $I,J \subset S$, let $R= G_Ix_0$ and $Q = gG_J x_0$. Then 
 \[ G_I\cap g G_J g^{-1} = \gamma G_K \gamma^{-1}, \]
 where $\gamma \in G_I$ and  $K=\{s\in I : w^{-1} s w = t \text{ for some } t \in J \}$ with  $w=w(R,Q)$.
\end{coro}
  
 \begin{proof}
 Let $P=  G_I\cap g G_J g^{-1}$,  we choose $\gamma \in G_I$  such that  $ \di{\gamma x_0,Q}=\di{R, Q}$. As in $\Delta$ parallel residues have equal stabilizers, with Proposition \ref{PropPrincip}   \[P=\mathrm{Stab}_G(R)\cap\mathrm{Stab}_G(Q) =\mathrm{Stab}_G(\proj{R}{Q})=\gamma G_K \gamma^{-1}.\]   
 On the other hand, the type $K$ of the residue $\proj{R}{Q}$ is given by Proposition \ref{propBaseResiduesPara}.$vi)$ which finishes the proof. 
 \end{proof}
 
  \bibliography{Biblio}

\begin{thebibliography}{MPW15}

\bibitem[AB08]{AbramBrown}
Peter Abramenko and Kenneth~S. Brown.
\newblock {\em Buildings}, volume 248 of {\em Graduate Texts in Mathematics}.
\newblock Springer, New York, 2008.
\newblock Theory and applications.

\bibitem[AM15]{AntolinMinasyanTitsAlterGP}
Yago Antol{\'{\i}}n and Ashot Minasyan.
\newblock Tits alternatives for graph products.
\newblock {\em J. Reine Angew. Math.}, 704:55--83, 2015.

\bibitem[Cap14]{CapraceAutomRightAngled}
Pierre-Emmanuel Caprace.
\newblock Automorphism groups of right-angled buildings: simplicity and local
  splittings.
\newblock {\em Fund. Math.}, 224(1):17--51, 2014.

\bibitem[Cha07]{CharneyRAAG}
Ruth Charney.
\newblock An introduction to right-angled {A}rtin groups.
\newblock {\em Geom. Dedicata}, 125:141--158, 2007.

\bibitem[Dav98]{DavisCAT}
Michael~W. Davis.
\newblock Buildings are {${\rm CAT}(0)$}.
\newblock In {\em Geometry and cohomology in group theory ({D}urham, 1994)},
  volume 252 of {\em London Math. Soc. Lecture Note Ser.}, pages 108--123.
  Cambridge Univ. Press, Cambridge, 1998.

\bibitem[Dav08]{DavisBook}
Michael~W. Davis.
\newblock {\em The geometry and topology of {C}oxeter groups}, volume~32 of
  {\em London Mathematical Society Monographs Series}.
\newblock Princeton University Press, Princeton, NJ, 2008.

\bibitem[GP01]{GabPauImmeubles}
Damien Gaboriau and Fr{\'e}d{\'e}ric Paulin.
\newblock Sur les immeubles hyperboliques.
\newblock {\em Geom. Dedicata}, 88(1-3):153--197, 2001.

\bibitem[HP03]{HaglundPaulinImmeubles}
Fr{\'e}d{\'e}ric Haglund and Fr{\'e}d{\'e}ric Paulin.
\newblock Constructions arborescentes d'immeubles.
\newblock {\em Math. Ann.}, 325(1):137--164, 2003.

\bibitem[MPW15]{MuhlherrPeterWeissDescentinBuildings}
Bernhard M{\"u}hlherr, Holger~P. Petersson, and Richard~M. Weiss.
\newblock {\em Descent in buildings}, volume 190 of {\em Annals of Mathematics
  Studies}.
\newblock Princeton University Press, Princeton, NJ, 2015.

\bibitem[Ron89]{RonanBuildings}
Mark Ronan.
\newblock {\em Lectures on buildings}, volume~7 of {\em Perspectives in
  Mathematics}.
\newblock Academic Press, Inc., Boston, MA, 1989.

\bibitem[Tit74]{TitsBuildingsLectureNotes}
Jacques Tits.
\newblock {\em Buildings of spherical type and finite {BN}-pairs}.
\newblock Lecture Notes in Mathematics, Vol. 386. Springer-Verlag, Berlin-New
  York, 1974.

\bibitem[Tit92]{TitsTwinBuildingsGroupsKacMoody}
Jacques Tits.
\newblock Twin buildings and groups of {K}ac-{M}oody type.
\newblock In {\em Groups, combinatorics \& geometry ({D}urham, 1990)}, volume
  165 of {\em London Math. Soc. Lecture Note Ser.}, pages 249--286. Cambridge
  Univ. Press, Cambridge, 1992.

\bibitem[Wei03]{WeissSphericalBuildings}
Richard~M. Weiss.
\newblock {\em The structure of spherical buildings}.
\newblock Princeton University Press, Princeton, NJ, 2003.

\end{thebibliography}
  \bibliographystyle{alpha}

\end{document}